\documentclass[11pt]{amsart}
\usepackage{amsmath,amssymb}

\newtheorem{theorem}{Theorem}
\newtheorem{proposition}[theorem]{Proposition}
\newtheorem{corollary}[theorem]{Corollary}

\begin{document}

\title{A note on Ricci flow and optimal transportation}
\author{Simon Brendle}
\address{Department of Mathematics \\ Stanford University \\ Stanford, CA 94305}
\thanks{The author was partially supported by the National Science Foundation under grant DMS-0905628.}
\maketitle 

\section{Introduction}

In this note, we describe an interpolation inequality in the setting of Ricci flow and $\mathcal{L}$-distance. This inequality is motivated by the following classical inequality due to Pr\'ekopa and Leindler: 

\begin{theorem}[Pr\'ekopa \cite{Prekopa}; Leindler \cite{Leindler}]
Fix a real number $0 < \lambda < 1$. Moreover, let $u_1,u_2,v: \mathbb{R}^n \to \mathbb{R}$ be nonnegative measurable functions satisfying 
\[v((1-\lambda)x+\lambda y) \geq u_1(x)^{1-\lambda} \, u_2(y)^\lambda\] 
for all points $x,y \in \mathbb{R}^n$. Then 
\[\int_{\mathbb{R}^n} v \geq \bigg ( \int_{\mathbb{R}^n} u_1 \bigg )^\lambda \, \bigg ( \int_{\mathbb{R}^n} u_2 \bigg )^{1-\lambda}.\]
\end{theorem}

Cordero-Erausquin, McCann, and Schmuckenschl\"ager \cite{Cordero-Erausquin-McCann-Schmuckenschlager1}, \cite{Cordero-Erausquin-McCann-Schmuckenschlager2} have generalized this inequality to Riemannian manifolds. The proof employs techniques from optimal transport theory (see e.g. \cite{Villani}).

In this note, we replace the Riemannian distance by Perelman's $\mathcal{L}$-distance (cf. \cite{Perelman}). The theory of $\mathcal{L}$-optimal transport was developed in recent work of Topping \cite{Topping2} (see also \cite{Lott}, \cite{McCann-Topping}). Among other things, Topping proved an important monotonicity formula for the $\mathcal{L}$-Wasserstein distance on the space of probability measures. Lott \cite{Lott} established a convexity property for the entropy along $\mathcal{L}$-Wasserstein geodesics.

To fix notation, let $M$ be a compact manifold of dimension $n$, and let $g(t)$, $t \in [0,T]$, be a one-parameter family of metrics on $M$. We assume that the metrics $g(t)$ evolve by backward Ricci flow, i.e. 
\[\frac{\partial}{\partial t} g(t) = 2 \, \text{\rm Ric}_{g(t)}.\] 
This evolution equation was introduced in a famous paper by R.~Hamilton \cite{Hamilton}. For an introduction to Ricci flow, see \cite{Brendle} or \cite{Topping1}. Following Perelman \cite{Perelman}, we define the $\mathcal{L}$-length of a path $\gamma: [\tau_1,\tau_2] \to M$ by 
\[\mathcal{L}(\gamma) = \int_{\tau_1}^{\tau_2} \sqrt{t} \, (\text{\rm scal}_{g(t)}(\gamma(t)) + |\gamma'(t)|_{g(t)}^2) \, dt.\] 
Moreover, the $\mathcal{L}$-distance is defined by 
\[Q(x,\tau_1;y,\tau_2) = \inf \{\mathcal{L}(\gamma): \gamma: [\tau_1,\tau_2] \to M, \, \gamma(\tau_1)=x, \, \gamma(\tau_2)=y\}.\] 
A path $\gamma: [\tau_1,\tau_2] \to M$ is called an $\mathcal{L}$-geodesic if the first variation of $\mathcal{L}$ is zero. For each tangent vector $Z \in T_x M$, we define 
\[\mathcal{L}_{\tau_1,\tau_2} \exp_x(Z) = \gamma(\tau_2),\] 
where $\gamma: [\tau_1,\tau_2] \to M$ is the unique $\mathcal{L}$-geodesic satisfying $\gamma(\tau_1) = x$ and $\sqrt{\tau_1} \, \gamma'(\tau_1) = Z$. The following is the main result of this note: 

\begin{theorem}
\label{main.theorem}
Fix real numbers $\tau_1,\tau_2,\tau$ such that $0 < \tau_1 < \tau < \tau_2 < T$. For abbreviation, we write 
\[\frac{1}{\sqrt{\tau}} = \frac{1-\lambda}{\sqrt{\tau_1}} + \frac{\lambda}{\sqrt{\tau_2}},\] 
where $0 < \lambda < 1$. Let $u_1,u_2,v: M \to \mathbb{R}$ be nonnegative measurable functions such that 
\begin{align*} 
\Big ( \frac{\tau}{\tau_1^{1-\lambda} \, \tau_2^\lambda} \Big )^{\frac{n}{2}} \, v(\gamma(\tau)) 
&\geq \exp \Big ( -\frac{1-\lambda}{2\sqrt{\tau_1}} \, Q(\gamma(\tau_1),\tau_1;\gamma(\tau),\tau) \Big ) \, u_1(\gamma(\tau_1))^{1-\lambda} \\ 
&\cdot \exp \Big ( \frac{\lambda}{2\sqrt{\tau_2}} \, Q(\gamma(\tau),\tau;\gamma(\tau_2),\tau_2) \Big ) \, u_2(\gamma(\tau_2))^\lambda 
\end{align*} 
for every minimizing $\mathcal{L}$-geodesic $\gamma: [\tau_1,\tau_2] \to M$. Then 
\[\int_M v \, d\text{\rm vol}_{g(\tau)} \geq \bigg ( \int_M u_1 \, d\text{\rm vol}_{g(\tau_1)} \bigg )^{1-\lambda} \, \bigg ( \int_M u_2 \, d\text{\rm vol}_{g(\tau_2)} \bigg )^\lambda.\] 
\end{theorem}

By sending $\tau_1 \to 0$, we recover the monotonicity of Perelman's reduced volume. This is discussed in Section \ref{relation.to.reduced.volume}.

\section{Proof of Theorem \ref{main.theorem}}

In order to prove Theorem \ref{main.theorem}, we make extensive use of Topping's notion of $\mathcal{L}$-optimal transportation (cf. \cite{Topping2}, \cite{Lott}). Without loss of generality, we may assume that
\[\int_M u_1(x) \, d\text{\rm vol}_{g(\tau_1)} = \int_M u_2(x) \, d\text{\rm vol}_{g(\tau_2)} = 1.\]
We define Borel probability measures $\nu_1$ and $\nu_2$ by
\[d\nu_1 = u_1 \, d\text{\rm vol}_{g(\tau_1)}\]
and
\[d\nu_2 = u_2 \, d\text{\rm vol}_{g(\tau_2)}.\] 
The following result was proved by Topping (cf. \cite{Topping2}, Section 2; see also \cite{Villani}, Theorem 10.28).

\begin{proposition}[P.~Topping \cite{Topping2}]
\label{existence.of.optimal.transport}
There exists a Borel map $F: M \to M$, a reflexive function $\varphi: M \to \mathbb{R}$, and a Borel set $K \subset M$ with the following properties: 
\begin{itemize}
\item[(i)] The set $M \setminus K$ has measure zero.
\item[(ii)] If $x \in K$, then the functions $\varphi$ and $Q(\cdot,\tau_1;F(x),\tau_2)$ are both differentiable at $x$. Moreover, the gradient of the function $Q(\cdot,\tau_1;F(x),\tau_2) - \varphi$ at the point $x$ is equal to zero.
\item[(iii)] We have $\nu_2 = F_\# \nu_1$.
\end{itemize}
\end{proposition}

Since $\varphi$ is reflexive, the function $\varphi$ is semiconcave (cf. \cite{Topping2}, Lemma 2.10). By Theorem 14.1 in \cite{Villani}, we can find a Borel set $\tilde{K} \subset K$ with the following properties:
\begin{itemize}
\item The set $M \setminus \tilde{K}$ has measure zero.
\item For each point $x \in \tilde{K}$, the function $\varphi$ admits a Taylor expansion of order two around $x$.
\end{itemize}
For each point $x \in \tilde{K}$, we denote by $\nabla \varphi(x)$ and $(\text{\rm Hess} \, \varphi)_x$ the gradient and Hessian of the function $\varphi$ with respect to the metric $g(\tau_1)$. Theorem 14.1 in \cite{Villani} guarantees that $(\text{\rm Hess} \, \varphi)_x$ is symmetric.

For each $t \in [\tau_1,\tau_2]$, we define a Borel map $F_t: M \to M$ by 
\[F_t(x) = \mathcal{L}_{\tau_1,t} \exp_x \Big ( -\frac{1}{2} \, \nabla \varphi(x) \Big )\] 
for $x \in \tilde{K}$. The following result is a consequence of property (ii) in Proposition \ref{existence.of.optimal.transport}:

\begin{proposition}
\label{L.geodesic}
We have $F_{\tau_1}(x) = x$ and $F_{\tau_2}(x) = F(x)$ for all $x \in \tilde{K}$. Moreover, for each point $x \in \tilde{K}$, the path $t \mapsto F_t(x)$ has minimal $\mathcal{L}$-length among all paths joining $(x,\tau_1)$ and $(F(x),\tau_2)$.
\end{proposition}

We next describe the volume distortion coefficients. To that end, we fix a point $x \in \tilde{K}$ and time $t \in (\tau_1,\tau_2]$. The linearization of the $\mathcal{L}$-exponential map $\mathcal{L}_{\tau_1,t} \exp_x$ gives a linear transformation 
\[D(\mathcal{L}_{\tau_1,t} \exp_x)_{-\frac{1}{2} \, \nabla \varphi(x)}:  (T_x M,g(\tau_1)) \to (T_{F_t(x)} M,g(t)).\] 
Moreover, the Hessian of the function $Q(\cdot,\tau_1;F_t(x),t) - \varphi$ at the point $x$ defines a symmetric linear transformation from the tangent space $(T_x M,g(\tau_1))$ into itself. 

Let us denote by $\Psi_{x,t}: (T_x M,g(\tau_1)) \to (T_{F_t(x)} M,g(t))$ the composition of these two linear transformations; that is, 
\[\Psi_{x,t} = \frac{1}{2} \, D(\mathcal{L}_{\tau_1,t} \exp_x)_{-\frac{1}{2} \, \nabla \varphi(x)} \circ \Big [ \text{\rm Hess}(Q(\cdot,\tau_1;F_t(x),t) - \varphi) \Big ]_x\] 
(cf. \cite{Topping2}, Lemma 2.13). We then define the volume distortion coefficients by 
\[\mathcal{J}(x,t) = \det \Psi_{x,t}\] 
for all $x \in \tilde{K}$ and all $t \in (\tau_1,\tau_2]$.

\begin{proposition}
\label{Jacobian}
For each point $x \in \tilde{K}$, we have 
\begin{align*} 
&\tau^{-\frac{n}{2}} \, \exp \Big ( -\frac{1-\lambda}{2\sqrt{\tau_1}} \, Q(x,\tau_1;F_\tau(x),\tau) \Big ) \, \mathcal{J}(x,\tau) \\ 
&\geq \tau_1^{-\frac{n(1-\lambda)}{2}} \, \tau_2^{-\frac{n\lambda}{2}} \, \exp \Big ( -\frac{\lambda}{2\sqrt{\tau_2}} \, Q(F_\tau(x),\tau;F(x),\tau_2) \Big ) \, \mathcal{J}(x,\tau_2)^\lambda. 
\end{align*}
\end{proposition}

\textbf{Proof.}
Fix a point $x \in \tilde{K}$, and let $\{\bar{e}_1,\hdots,\bar{e}_n\}$ be a basis of $T_x M$ which is orthonormal with respect to the metric $g(\tau_1)$. We next consider the path 
\[\gamma: [\tau_1,\tau_2] \to M, \quad t \mapsto F_t(x).\] 
Let $\{e_1(t), \hdots, e_n(t)\}$ be vector fields along $\gamma$ such that $e_i(\tau_1) = \overline{e}_i$ and 
\[\langle D_t e_i(t),e_j(t) \rangle_{g(t)} + \text{\rm Ric}_{g(t)}(e_i(t),e_j(t)) = 0\] 
for all $t \in [\tau_1,\tau_2]$. For each $t \in [\tau_1,\tau_2]$, the vectors $\{e_1(t),\hdots,e_n(t)\}$ are orthonormal with respect to the metric $g(t)$.

Let $\{Y_1(t),\hdots,Y_n(t)\}$ be $\mathcal{L}$-Jacobi fields along $\gamma$ such that $Y_j(\tau_1) = \bar{e}_j$ and 
\[\langle \bar{e}_i,D_t Y_j(\tau_1) \rangle_{g(\tau_1)} = -\frac{1}{2\sqrt{\tau_1}} \, (\text{\rm Hess} \, \varphi)_x(\bar{e}_i,\bar{e}_j).\] 
For each $t \in [\tau_1,\tau_2]$, we define an $n \times n$ matrix $A(t)$ by
\[a_{ij}(t) = \langle e_i(t),Y_j(t) \rangle_{g(t)}.\] 
It follows from the initial conditions for $Y_j$ that $a_{ij}(\tau_1) = \delta_{ij}$ and 
\[a_{ij}'(\tau_1) = \text{\rm Ric}_{g(\tau_1)}(\bar{e}_i,\bar{e}_j) - \frac{1}{2\sqrt{\tau_1}} \, (\text{\rm Hess} \, \varphi)_x(\bar{e}_i,\bar{e}_j).\] 
In particular, the matrix $A'(\tau_1) \, A(\tau_1)^{-1}$ is symmetric. Moreover, it was shown by Topping \cite{Topping2} that 
\[A''(t) + \frac{1}{2t} \, A'(t) = M(t) A(t)\]
for all $t \in [\tau_1,\tau_2]$. Here, $M(t)$ is a symmetric $n \times n$ matrix, whose trace is given by 
\begin{align*} 
2 \, \text{\rm tr}(M(t)) 
&= \frac{\partial}{\partial t} \text{\rm scal}_{g(t)}(\gamma(t)) + 2 \, \langle \nabla \text{\rm scal}_{g(t)}(\gamma(t)),\gamma'(t) \rangle_{g(t)} \\ 
&- 2 \, \text{\rm Ric}_{g(t)}(\gamma'(t),\gamma'(t)) + \frac{1}{t} \, \text{\rm scal}_{g(t)}(\gamma(t)). 
\end{align*} 
Since $M(t)$ is symmetric, we conclude that the matrix $A'(t) \, A(t)^{-1}$ is symmetric for all $t \in [\tau_1,\tau_2]$. Arguing as in the proof of Lemma 3.1 in \cite{Topping2}, we obtain 
\begin{align*}
&t^{-\frac{3}{2}} \, \frac{d}{dt} \Big [ t^{\frac{3}{2}} \, \frac{d}{dt} \log \det A(t) \Big ] \\ 
&= \frac{d^2}{dt^2} \log \det A(t) + \frac{3}{2t} \, \frac{d}{dt} \log \det A(t) \\ 
&= \text{\rm tr}(A''(t) \, A(t)^{-1}) - \text{\rm tr}(A'(t) \, A(t)^{-1} \, A'(t) \, A(t)^{-1}) + \frac{3}{2t} \, \text{\rm tr}(A'(t) \, A(t)^{-1}) \\ 
&= \text{\rm tr}(M(t)) - \text{\rm tr}(A'(t) \, A(t)^{-1} \, A'(t) \, A(t)^{-1}) + \frac{1}{t} \, \text{\rm tr}(A'(t) \, A(t)^{-1}) \\ 
&= \text{\rm tr}(M(t)) - \text{\rm tr} \Big [ \big ( A'(t) \, A(t)^{-1} - \frac{1}{2t} \, I \big )^2 \Big ] + \frac{n}{4t^2} \\ 
&\leq \text{\rm tr}(M(t)) + \frac{n}{4t^2} 
\end{align*}
for all $t \in [\tau_1,\tau_2]$. Moreover, we have 
\[\frac{d}{dt} Q(x,\tau_1;F_t(x),t) = \sqrt{t} \, (\text{\rm scal}_{g(t)}(\gamma(t)) + |\gamma'(t)|_{g(t)}^2)\] 
by definition of the $\mathcal{L}$-distance. This implies 
\begin{align*} 
&t^{-\frac{3}{2}} \, \frac{d}{dt} \Big [ t^{\frac{3}{2}} \, \frac{d}{dt} \big ( t^{-\frac{1}{2}} \, Q(x,\tau_1;F_t(x),t) \big ) \Big ] \\ 
&= t^{-1} \, \frac{d}{dt} \Big [ t^{\frac{1}{2}} \, \frac{d}{dt} Q(x,\tau_1;F_t(x),t) \Big ] \\ 
&= t^{-1} \, \frac{d}{dt} \Big [ t \, (\text{\rm scal}_{g(t)}(\gamma(t)) + |\gamma'(t)|_{g(t)}^2) \Big ] \\ 
&= \frac{\partial}{\partial t} \text{\rm scal}_{g(t)}(\gamma(t)) + 2 \, \langle \nabla \text{\rm scal}_{g(t)}(\gamma(t)),\gamma'(t) \rangle_{g(t)} \\ 
&- 2 \, \text{\rm Ric}_{g(t)}(\gamma'(t),\gamma'(t)) + \frac{1}{t} \, \text{\rm scal}_{g(t)}(\gamma(t)) \\ 
&= 2 \, \text{\rm tr}(M(t)) 
\end{align*}
(cf. \cite{Perelman}, equation (7.3)). Putting these facts together, we obtain 
\[t^{-\frac{3}{2}} \, \frac{d}{dt} \Big [ t^{\frac{3}{2}} \, \frac{d}{dt} \Big ( \frac{n}{2} \, \log t + \frac{1}{2} \, t^{-\frac{1}{2}} \, Q(x,\tau_1;F_t(x),t) - \log \det A(t) \Big ) \Big ] \geq 0.\] 
Hence, if we write 
\[\frac{n}{2} \, \log t + \frac{1}{2} \, t^{-\frac{1}{2}} \, Q(x,\tau_1;F_t(x),t) - \log \det A(t) = h(t^{-\frac{1}{2}}),\] 
then the function $h$ is convex. Since $\tau^{-\frac{1}{2}} = (1-\lambda) \, \tau_1^{-\frac{1}{2}} + \lambda \, \tau_2^{-\frac{1}{2}}$, we conclude that  
\[h(\tau^{-\frac{1}{2}}) \leq (1-\lambda) \, h(\tau_1^{-\frac{1}{2}}) + \lambda \, h(\tau_2^{-\frac{1}{2}}),\]
hence 
\begin{align*} 
&\tau^{-\frac{n}{2}} \, \exp \Big ( -\frac{1}{2\sqrt{\tau}} \, Q(x,\tau_1;F_\tau(x),\tau) \Big ) \, \det A(\tau) \\ 
&\geq \tau_1^{-\frac{n(1-\lambda)}{2}} \, \tau_2^{-\frac{n\lambda}{2}} \, \exp \Big ( -\frac{\lambda}{2\sqrt{\tau_2}} \, Q(x,\tau_1;F(x),\tau_2) \Big ) \, (\det A(\tau_2))^\lambda. 
\end{align*}
By Lemma \ref{L.geodesic}, the path $\gamma$ has minimal $\mathcal{L}$-length among all paths joining $(x,\tau_1)$ to $(F(x),\tau_2)$. This implies 
\[Q(x,\tau_1;F(x),\tau_2) = Q(x,\tau_1;F_\tau(x),\tau) + Q(F_\tau(x),\tau;F(x),\tau_2).\] 
Hence, we obtain 
\begin{align*} 
&\tau^{-\frac{n}{2}} \, \exp \Big ( -\frac{1-\lambda}{2\sqrt{\tau_1}} \, Q(x,\tau_1;F_\tau(x),\tau) \Big ) \, \det A(\tau) \\ 
&\geq \tau_1^{-\frac{n(1-\lambda)}{2}} \, \tau_2^{-\frac{n\lambda}{2}} \, \exp \Big ( -\frac{\lambda}{2\sqrt{\tau_2}} \, Q(F_\tau(x),\tau;F(x),\tau_2) \Big ) \, (\det A(\tau_2))^\lambda. 
\end{align*}
On the other hand, it follows from Lemma 2.18 in \cite{Topping2} that $\Psi_{x,t}(\bar{e}_j) = Y_j(t)$ for all $t \in (\tau_1,\tau_2]$. From this, we deduce that $\langle e_i(t),\Psi_{x,t}(\bar{e}_j) \rangle_{g(t)} = a_{ij}(t)$, hence $\mathcal{J}(x,t) = \det A(t)$ for all $t \in (\tau_1,\tau_2]$. Putting these facts together, the assertion follows. This completes the proof of Proposition \ref{Jacobian}. \\

We next consider the measure $\nu = (F_\tau)_\# \nu_1$. It follows from general results in optimal transportation theory (see e.g. \cite{Villani}) that $\nu$ is absolutely continuous with respect to the volume measure, i.e. $d\nu = u \, d\text{\rm vol}_{g(\tau)}$.

\begin{proposition}
\label{interpolant.density}
There exists a Borel set $\hat{K} \subset \tilde{K}$ such that $M \setminus \hat{K}$ has measure zero and
\begin{align*}
\Big ( \frac{\tau}{\tau_1^{1-\lambda} \, \tau_2^\lambda} \Big )^{\frac{n}{2}} \, u(F_\tau(x)) 
&\leq \exp \Big ( -\frac{1-\lambda}{2\sqrt{\tau_1}} \, Q(x,\tau_1;F_\tau(x),\tau) \Big ) \, u_1(x)^{1-\lambda} \\ 
&\cdot \exp \Big ( \frac{\lambda}{2\sqrt{\tau_2}} \, Q(F_\tau(x),\tau;F(x),\tau_2) \Big ) \, u_2(F(x))^\lambda 
\end{align*} 
for all $x \in \hat{K}$.
\end{proposition}

\textbf{Proof.}
It follows from Theorem 2.14 in \cite{Topping2} that
\[u_1(x) = u_2(F(x)) \, \mathcal{J}(x,\tau_2) > 0\]
for almost all $x \in \tilde{K}$. An analogous argument yields
\[u_1(x) = u(F_\tau(x)) \, \mathcal{J}(x,\tau) > 0\]
for almost all $x \in \tilde{K}$. Using Proposition \ref{Jacobian}, we obtain 
\begin{align*} 
&\tau^{-\frac{n}{2}} \, \exp \Big ( -\frac{1-\lambda}{2\sqrt{\tau_1}} \, Q(x,\tau_1;F_\tau(x),\tau) \Big ) \, \frac{u_1(x)}{u(F_\tau(x))} \\ 
&\geq \tau_1^{-\frac{n(1-\lambda)}{2}} \, \tau_2^{-\frac{n\lambda}{2}} \, \exp \Big ( -\frac{\lambda}{2\sqrt{\tau_2}} \, Q(F_\tau(x),\tau;F(x),\tau_2) \Big ) \, \Big ( \frac{u_1(x)}{u_2(F(x))} \Big )^\lambda 
\end{align*}
for almost all $x \in \tilde{K}$. Rearranging terms, the assertion follows. \\

\begin{corollary}
We have 
\[\int_M v \, d\text{\rm vol}_{g(\tau)} \geq 1.\] 
\end{corollary}

\textbf{Proof.}
Fix a point $x \in \hat{K}$. By Lemma \ref{L.geodesic}, the path $t \mapsto F_t(x)$ is a minimizing $\mathcal{L}$-geodesic. Therefore, we have 
\begin{align*} 
\Big ( \frac{\tau}{\tau_1^{1-\lambda} \, \tau_2^\lambda} \Big )^{\frac{n}{2}} \, v(F_\tau(x)) 
&\geq \exp \Big ( -\frac{1-\lambda}{2\sqrt{\tau_1}} \, Q(x,\tau_1;F_\tau(x),\tau) \Big ) \, u_1(x)^{1-\lambda} \\ 
&\cdot \exp \Big ( \frac{\lambda}{2\sqrt{\tau_2}} \, Q(F_\tau(x),\tau;F(x),\tau_2) \Big ) \, u_2(F(x))^\lambda. 
\end{align*}
Using Proposition \ref{interpolant.density}, we conclude that 
\[v(F_\tau(x)) \geq u(F_\tau(x))\]
for all $x \in \hat{K}$. This implies
\begin{align*}
\int_M v \, d\text{\rm vol}_{g(\tau)}
&\geq \int_{F_\tau(\hat{K})} v \, d\text{\rm vol}_{g(\tau)} \\
&\geq \int_{F_\tau(\hat{K})} u \, d\text{\rm vol}_{g(\tau)} \\
&= \nu(F_\tau(\hat{K})).
\end{align*}
Moreover, we have
\[\nu(F_\tau(\hat{K})) = \nu_1 \big [ F_\tau^{-1}(F_\tau(\hat{K})) \big ] \geq \nu_1(\hat{K}) = 1\]
by definition of $\nu$. Putting these facts together, the assertion follows. \\

\section{Relation to Perelman's reduced volume}

\label{relation.to.reduced.volume}

In this final section, we discuss how Theorem \ref{main.theorem} is related to the monotonicity of Perelman's reduced volume. The strategy is to fix $\tau$ and $\tau_2$, and pass to the limit as $\tau_1 \to 0$.

Let us fix a point $p \in M$ and real numbers $0 < \tau < \tau_2 < T$. We define a function $v$ by 
\[v = \tau^{-\frac{n}{2}} \, \exp \Big ( -\frac{1}{2\sqrt{\tau}} \, Q(p,0;\cdot,\tau) \Big ).\] 
For $\tau_1 > 0$ sufficiently small, we denote by $B(p,\sqrt{\tau_1})$ the geodesic ball of radius $\sqrt{\tau_1}$ in the metric $g(0)$. We can find a positive constant $N$ such that $Q(p,0;x,\tau_1) \leq N \, \sqrt{\tau_1}$ and $Q(x,\tau_1;p,2\tau_1) \leq N \, \sqrt{\tau_1}$ for all points $x \in B(p,\sqrt{\tau_1})$. Note that the constant $N$ is independent of $\tau_1$.

As above, we write 
\[\frac{1}{\sqrt{\tau}} = \frac{1-\lambda}{\sqrt{\tau_1}} + \frac{\lambda}{\sqrt{\tau_2}},\] 
where $0 < \lambda < 1$. We now specify the functions $u_1$ and $u_2$. We define 
\[u_1 = \tau_1^{-\frac{n}{2}} \, \exp \Big ( -\frac{N\sqrt{\tau_1}}{2(1-\lambda)} \, \Big ( \frac{1}{\sqrt{\tau}} + \frac{\lambda}{\sqrt{\tau_2}} \Big ) \Big ) \, 1_{B(p,\sqrt{\tau_1})}\] 
and 
\[u_2 = \tau_2^{-\frac{n}{2}} \, \exp \Big ( -\frac{1}{2\sqrt{\tau_2}} \, Q(p,2\tau_1;\cdot,\tau_2) \Big ).\] 
We claim that $u_1,u_2,v$ satisfy the assumptions of Theorem \ref{main.theorem}:

\begin{proposition}
We have 
\begin{align*} 
\Big ( \frac{\tau}{\tau_1^{1-\lambda} \, \tau_2^\lambda} \Big )^{\frac{n}{2}} \, v(\gamma(\tau)) 
&\geq \exp \Big ( -\frac{1-\lambda}{2\sqrt{\tau_1}} \, Q(\gamma(\tau_1),\tau_1;\gamma(\tau),\tau) \Big ) \, u_1(\gamma(\tau_1))^{1-\lambda} \\ 
&\cdot \exp \Big ( \frac{\lambda}{2\sqrt{\tau_2}} \, Q(\gamma(\tau),\tau;\gamma(\tau_2),\tau_2) \Big ) \, u_2(\gamma(\tau_2))^\lambda 
\end{align*} 
for every minimizing $\mathcal{L}$-geodesic $\gamma: [\tau_1,\tau_2] \to M$.
\end{proposition}

\textbf{Proof.} 
If $\gamma(0) \notin B(p,\sqrt{\tau_1})$, then $u_1(\gamma(0)) = 0$ and the assertion is trivial. Hence, it suffices to consider the case $\gamma(0) \in B(p,\sqrt{\tau_1})$. In this case, we have 
\begin{align*} 
Q(p,0;\gamma(\tau),\tau) 
&\leq Q(\gamma(\tau_1),\tau_1;\gamma(\tau),\tau) + Q(p,0;\gamma(\tau_1),\tau_1) \\ 
&\leq Q(\gamma(\tau_1),\tau_1;\gamma(\tau),\tau) + N \, \sqrt{\tau_1} 
\end{align*} 
and 
\begin{align*} 
Q(p,2\tau_1;\gamma(\tau_2),\tau_2) 
&\geq Q(\gamma(\tau_1),\tau_1;\gamma(\tau_2),\tau_2) - Q(\gamma(\tau_1),\tau_1;p,2\tau_1) \\ 
&\geq Q(\gamma(\tau_1),\tau_1;\gamma(\tau_2),\tau_2) - N \, \sqrt{\tau_1}. 
\end{align*} 
This implies 
\[v(\gamma(\tau),\tau) \geq \tau^{-\frac{n}{2}} \, \exp \Big ( -\frac{N \, \sqrt{\tau_1}}{2\sqrt{\tau}} \Big ) \, \exp \Big ( -\frac{1}{2\sqrt{\tau}} \, Q(\gamma(\tau_1),\tau_1;\gamma(\tau),\tau) \Big )\] 
and 
\[u_2(\gamma(\tau_2)) \leq \tau_2^{-\frac{n}{2}} \, \exp \Big ( \frac{N \, \sqrt{\tau_1}}{2\sqrt{\tau_2}} \Big ) \, \exp \Big ( -\frac{1}{2\sqrt{\tau_2}} \, Q(\gamma(\tau_1),\tau_1;\gamma(\tau_2),\tau_2) \Big ).\] 
Moreover, we have
\[Q(\gamma(\tau_1),\tau_1;\gamma(\tau_2),\tau_2) = Q(\gamma(\tau_1),\tau_1;\gamma(\tau),\tau) + Q(\gamma(\tau),\tau;\gamma(\tau_2),\tau_2)\] 
since $\gamma$ has minimal $\mathcal{L}$-length. Putting these facts together, we obtain 
\begin{align*} 
v(\gamma(\tau)) 
&\geq \Big ( \frac{\tau_2^\lambda}{\tau} \Big )^{\frac{n}{2}} \, \exp \Big ( -\frac{N\sqrt{\tau_1}}{2} \, \Big ( \frac{1}{\sqrt{\tau}} + \frac{\lambda}{\sqrt{\tau_2}} \Big ) \Big ) \\ 
&\cdot \exp \Big ( -\frac{1-\lambda}{2\sqrt{\tau_1}} \, Q(\gamma(\tau_1),\tau_1;\gamma(\tau),\tau) \Big ) \\ 
&\cdot \exp \Big ( \frac{\lambda}{2\sqrt{\tau_2}} \, Q(\gamma(\tau),\tau;\gamma(\tau_2),\tau_2) \Big ) \, u_2(\gamma(\tau_2))^\lambda. 
\end{align*} From this, the assertion follows. \\

Let $\tilde{V}(\tau)$ denote the reduced volume at time $\tau$. Using Theorem \ref{main.theorem}, we obtain 
\begin{align*} 
\tilde{V}(\tau) &= \int_M v \, d\text{\rm vol}_{g(\tau)} \\ 
&\geq \bigg ( \int_M u_1 \, d\text{\rm vol}_{g(\tau_1)} \bigg )^{1-\lambda} \, \bigg ( \int_M u_2 \, d\text{\rm vol}_{g(\tau_2)} \bigg )^\lambda. 
\end{align*}
We now fix $\tau$ and $\tau_2$, and pass to the limit as $\tau_1 \to 0$. Clearly, 
\[1-\lambda = \sqrt{\tau_1} \, \Big ( \frac{1}{\sqrt{\tau}} - \frac{1}{\sqrt{\tau_2}} \Big ) + O(\tau_1).\] 
This implies 
\[-\frac{N\sqrt{\tau_1}}{2(1-\lambda)} \, \Big ( \frac{1}{\sqrt{\tau}} + \frac{\lambda}{\sqrt{\tau_2}} \Big ) \to -\frac{N}{2} \, \Big ( \frac{1}{\sqrt{\tau}} - \frac{1}{\sqrt{\tau_2}} \Big )^{-1} \, \Big ( \frac{1}{\sqrt{\tau}} + \frac{1}{\sqrt{\tau_2}} \Big ).\] 
Hence, the integral $\int_M u_1 \, d\text{\rm vol}_{g(\tau_1)}$ converges to a positive real number as $\tau_1 \to 0$. Since $1-\lambda \to 0$, we conclude that 
\[\bigg ( \int_M u_1 \, d\text{\rm vol}_{g(\tau_1)} \bigg )^{1-\lambda} \to 1\] 
as $\tau_1 \to 0$. Moreover, we have 
\[\bigg ( \int_M u_2 \, d\text{\rm vol}_{g(\tau_2)} \bigg )^\lambda \to \tilde{V}(\tau_2)\] 
as $\tau_1 \to 0$. Putting these facts together, we obtain 
\[\tilde{V}(\tau) = \int_M v \, d\text{\rm vol}_{g(\tau)} \geq \tilde{V}(\tau_2).\] 
Thus, Theorem \ref{main.theorem} implies the monotonicity of the reduced volume.

\end{document}